\documentstyle[12pt]{article}
\begin{document}
\newtheorem{theorem}{Theorem}
\newtheorem{opr}{Definition}
\newtheorem{pr}{Propostion}
\newtheorem{remark}{Remark}
\newtheorem{question}{Question}
\newtheorem{lemma}{Lemma}
\newtheorem{co}{Corollary}
\newtheorem{example}{Example}

\begin{center}
\Large{
{\bf Exponential groups 2:\\
Extensions of centralizers and tensor completion of CSA--groups}\\
by \\
Myasnikov A.G., Remeslennikov V.N.}
\end{center}

This is the second article in the series of papers by the authors
on the theory of exponential groups. In the first one \cite{mr1} we 
discussed foundations of this theory. 
Definitions necessary for independent understanding of the present article
are given in the introduction and the first section.
The theory of exponential groups 
begins with  results of A.Mal'cev \cite{mal}, P.Hall \cite{hall},
G.Baumslag \cite{bau} and  R.Lyndon \cite{lyn}.
The  axiomatic notion of an exponential group was introduced by R.Lyndon (1960).
In \cite{myas} a new axiom was added to  Lyndon's definition to obtain a new notion of
an exponential group. The refined version is more convenient because it
coincides exactly with the notion of a module over a ring in the
abelian case, whereas abelian exponential groups after Lyndon provide
a far wider class.\par
Recall the main definition from \cite{myas}.\par
Let $A$ be an arbitrary associative ring with identity and $G$ a group.
Fix an action of the ring $A$ on $G$, i.e. a map $G \times A \rightarrow
G$. The result of the action of $\alpha \in A$ on $g \in G$ is written
as $g^\alpha$. Consider the following axioms:

\begin{enumerate}
\item $g^1=g$, $g^0=1$,  $1^\alpha = 1$ ;
\item $g^{\alpha +\beta}=g^\alpha \cdot g^\beta, \ g^{\alpha \beta} =
(g^\alpha)^\beta$;
\item $(h^{-1}gh)^\alpha = h^{-1}g^\alpha h;$
\item $[g,h]=1 \Longrightarrow (gh)^\alpha = g^\alpha h^\alpha.$
\end{enumerate}

\begin{opr}
Groups with axioms 1)--4) are called $A$--groups. 
\end{opr}

As usual, one can introduce the notions of $A$--torsion, finitely 
$A$-generated group, $A$--homomorphism, $A$--isomorphism, and so on.\par

Our interest to the theory of exponential groups is motivated by the
following circumstances. First, many natural classes of groups are
$A$--groups. For example, unipotent groups over a field $k$ of
zero characteristic are $k$--groups, pro-$p$-groups are exponential groups 
over the ring of $p$--adic integers, etc. (see \cite{mr1} for examples).
Second, the notion of an $A$--group is a natural generalization of the
notion of a module over a ring in the category of noncommutative
groups. Third, this notions allows one to define for an arbitrary group
$G$ its largest ring of scalars $A(G)$, over which $G$ is an $A(G)$--group.
The last notion is an exact analog of the notion of the centroid for 
rings and algebras (see \cite{myas} for details) and plays an important part in 
many questions. Finally, the notion of an $A$--exponential group is 
helpful in investigation of model--theoretic problems for noncommutative
groups \cite{myas2}. Let us illustrate this on example of universally free
group, i.e. groups that has the same universal theory as a nonabelian free
group (see \cite{vnr} for details).
The study of these groups plays a key part in creating model theory of
free groups. Here we formulate a principal problem for universally free groups.
\par\medskip
{\bf Conjecture}  {\it Every finitely generated universally free group 
is a subgroup of a free ${\bf Z}[x]$--group, where ${\bf Z}[x]$
is the ring of polynomials in an indeterminate $x$ with integer
coefficients.}
\par\medskip
This problem requires a detailed study of ${\bf Z}[x]$--free groups
and their subgroups. \par 
	The present article will be devoted to 
the general theory of tensor  $A$--completions of groups with an emphasis on 
$A$--free groups..\par
	G.Baumslag was the first who introduced the ${\bf Q}$--completion  
$F^{\bf Q}$ of 
a free group $F$ and, what is more important, described it using free products with amalgamation. It allowed him to establish some properties of the group 
$F^{\bf Q}$. We develope his approach and apply it for arbitrary rings $A$ of characterisric 0 and groups from a very wide class.

One of the key notions in the theory of $A$--groups is the notion of tensor 
$A$--completion.
The case where a group $G$ is a subgroup of its $A$--completion 
$G^A$ (i.e. $G$ is $A$--faithful) is of a particular interest. 
In section 2 we prove that the class ${\cal F}_A$ of all $A$--faithful
groups is universally axiomatizable and closed under subgroups,
direct products and direct limits, and, in the case when $A$ is an integral
domain of characteristic zero, contains all residually torsion--free nilpotent groups. 
For arbitrary groups it is difficult to give a 
good constructive description of their $A$--completion. In sections 5--7
we study the class of CSA--groups, for which a good and concrete
description of
tensor completion exists. The class of CSA--groups is quite wide, 
it is universally axiomatizable, closed under free products,
direct limits and extensions of centralizers (basic facts concerning
these constructions are given in sections 3--4); moreover it 
contains abelian, free, hyperbolic torsion--free groups and groups
acting freely on $\Lambda$--trees and universally free groups. 
For CSA--groups the construction 
of tensor completion of a group $G$ is obtained as an iterated 
tree extension of centralizers of the initial group $G$. In 
particular, any torsion--free CSA--group is faithful over any ring
$A$ with a torsion--free additive group $A^+$. In the last
section we apply the obtained results to the study of basic 
properties of $A$--free groups. In particular, 
canonical and reduced forms of elements in $A$--free groups
are introduced, and then commuting and conjugate
elements are described. We finish the paper by formulating some
open problems on this area.\par
The results of this article are partially contained in preprints 
\cite{amr,mr2}.

\section{Tensor completions and free A--groups}

The basic operation in the class of $A$--groups is the operation of tensor
completion. Here we give a particular case of this construction
(see \cite{mr1} for the general definition). Later on we always assume that
the ring $A$ and its subring $A_0$ have a common identity element.

\begin{opr}
Let $G$ be an $A_0$--group and $A_0$ a subring of $A$. 
Then an $A$--group $G^A$ is called  a tensor $A$--completion
of the group $G$ if $G^A$ satisfies the following universal
property:
\begin{enumerate}
\item there exists an $A_0$--homomorphism $\lambda: G \rightarrow G^A$
such that $\lambda(G)$ $A$--generates $G^A$, i.e. 
$
\langle \lambda(G) \rangle_A = G^A
$
\item for any $A$--group $H$ and an $A_0$--homomorphism $\varphi: G \rightarrow H$
there exists the unique $A$--homomorphism $\psi: G^A \rightarrow H$
such that $\lambda \circ \psi \varphi$.
\end{enumerate}
\end{opr}

If $G$ is an abelian $A_0$--group, then the group $G^A$ is also
abelian, i.e. it is an $A$--module. $G^A$ satisfies the same
universal property as the tensor product $G \otimes_{A_0} A$ of the
$A_0$--module $G$ and the ring $A$. Therefore
$
G^A \simeq G \otimes_{A_0} A
$
\begin{theorem}{\rm \cite{mr1}}
Let $G$ be an $A_0$--group and $A_0$ a subring of $A$.
Then there exists a tensor completion $G^A$, and it is unique up to 
an $A$--isomorphism.
\end{theorem}
Every group is a $\bf Z$--group, so one can consider $A$--completions
of arbitrary groups for each ring of characteristic zero, i.e.
${\bf Z} \leq A$. \par

Let us formulate the notion of a free $A$-group. Let $A$ be an
associative ring with identity, char$A=0$, and $X$ an arbitrary set.

\begin{opr}
An $A$-group $F_A(X)$ with the set of $A$-generators $X$ is 
said to be a free $A$-group with base $X$, if for every
$A$-group $G$ an arbitrary mapping $\phi_0:X \mapsto G$ can
be extended to an $A$-homomorphism $\phi: F_A(X) \mapsto G$.
\end{opr}

The set $X$ is called the {\em set of free $A$-generators} of 
$F_A(X)$. The cardinality of $X$ is called the {\em rank} of
the group $F_A(X)$.

\begin{theorem}
For every $X$ and $A$ there exists a free $A$-group $F_A(X)$;
moreover it is unique up to $A$-isomorphism and $F_A(X)$ is
the tensor $A$-completion $F(X)^A$ of the group $F(X).$
\end{theorem}

$\Box$ Let $F(X)$ be a free group in the ordinary sense
and let $\phi_0: X \mapsto G$ be an 
arbitrary mapping from $X$ to an $A$-group $G$.
Then $\phi_0$ induces a homomorphism $\phi_1:F(X) \mapsto G$
by a property of a free group. Moreover, the latter mapping 
induces an $A$-homomorphism $\phi : (F(X))^A \mapsto G$.
Conseqently, $(F(X))^A$ is a free $A$-group with base $X$.
Uniqueness follows from the uniqueness of the tensor 
completion. $\Box$

\section{Faithful A--completion}

Usually one can construct the $A$--completion of a group $G$
step by step, on each stage creating an extension of a group
by defining an action of $A$ on some elements, and then closing
the process. Intermediate groups obtained along this construction
are partial $A$--groups. Following \cite{mr1}, let us remind the definition 
of a partial $A$--group.

\begin{opr} A group $G$ is a partial $A$--group iff $g^\alpha$
is defined for some $g \in G, \ \alpha \in A$, and axioms 1)--4) hold wherever
the action is defined (this means, in particular, that if one side of
equalities 1)--4) is defined, then so is the other side,and they are
equal to each other).
\end{opr}

Let $G$ be a partial $A$--group. One can define an $A$--{\em completion of a partial} $A$--{\em group} $G$ in exactly the same way as in definition 1
(considering $\lambda$ and $\phi$ to be partial $A$--homomorphisms).
In particular, if $A_0$ is a 
subring of $A$, then any $A_0$--group is also a partial 
$A$--group. So, consideration of tensor completions of partial $A$--groups 
covers the case of tensor $A$--completion of $A_0$--groups from section 1.
\begin{opr}
A partial $A$--group $G$ is called faithful (over $A$),
if the canonical homomorphism of partial $A$--groups 
$\lambda: G \longrightarrow G^A$ is injective. 
\end{opr}

\begin{remark}
\label{1}
\begin{enumerate}
\item A partial $A$--group $G$ is faithful iff $G$ is $A$--embeddible
into some $A$--group $H$.
\item  A group $G$ is  faithful over $A$ iff $G$ is residually faithful over $A$.
\end{enumerate}
\end{remark}
In \cite{am} it was shown that the operation of tensor completion 
commutes with the operation of direct product and direct limit:
$$
(\oplus G_i)^A \simeq_A \oplus G_i^A, \ \ \ 
(\lim_{\longrightarrow}G_i)^A \simeq_A \lim_{\longrightarrow}G_i^A
$$
Using these statements and the remark above, one can obtain the
following propositions:
\begin{pr}
Let ${\cal F}_A$ be a class of all $A$--faithful partial
$A$--groups. Then:
\begin{enumerate}
\item ${\cal F}_A$ is closed under direct products;
\item ${\cal F}_A$ is closed under direct limits.
\end{enumerate}
\end{pr}
$\Box$ Let $\lambda_i:G_i \longrightarrow G_i^A$ be canonical 
homomorphisms. By the properties above, $(\oplus G_i)^A \simeq
\oplus G_i^A$, and if $\lambda = \oplus \lambda_i$ is the direct sum
of $\lambda_i$, then the canonical map $\lambda: \oplus G_i \longrightarrow
\oplus G_i^A$ is injective iff each $\lambda_i$ is injective. This
implies 1). \par
Let $G = \lim_{\longrightarrow}G_i$. Then, as mentioned above,
$G^A = \lim_{\longrightarrow}G_i^A$. If $\lambda_i: G_i \rightarrow G_i^A$
is the canonical embedding, then there exists a canonical 
embedding
$$
\lambda: \lim_{\longrightarrow}G_i \longrightarrow \lim_{\longrightarrow}
G_i^A = G^A
$$
So $G=\lim_{\longrightarrow}G_i$ is $A$--faithful. \hfill $\Box$
\par
\medskip
The following proposition shows that the property of being faithful 
is a local property.
\begin{pr}
Let $G$ be a partial $A$--group. Then $G$ is faithful iff 
every finitely generated subgroup of $G$ is faithful.
\end{pr}
$\Box$
Let $\{ G_i, \ i \in I\}$ be the set of all finitely generated
partial $A$--subgroups of $G$. Then 
$G=\lim_{\longrightarrow}G_i$, and the statement follows from 
proposition 1. \hfill $\Box$

\begin{pr}
Let ${\cal F}_A$ be the class of all faithful partial $A$--groups.
Then:
\begin{enumerate}
\item ${\cal F}_A$ is closed under the operation of taking subgroups;
\item ${\cal F}_A$ is closed under ultraproducts;
\item ${\cal F}_A$ is universally axiomatizible.
\end{enumerate}
\end{pr}

$\Box$ 1) is evident by the remark above. To prove 2), one can consider
instead of the group $G$ the many-sorted structure 
${\cal A}(G,A)= \langle G, G^A, A \rangle$ with the group operations
on $G$ and $G^A$, ring operations on $A$, operations of action of 
$A$ on $G$ (partial) and on $G^A$, and the canonical map
$\lambda : G \longrightarrow G^A$. In ${\cal A}(G, A)$ the notion
of faithfulness is evidently axiomatizible by one sentence:
$g \neq 1 \Longrightarrow \lambda(g) \neq 1$. Instead of the
ultraproduct of groups $\Pi G_i /D$ over an ultrafilter $D$
one can consider the corresponding ultraproduct $\Pi {\cal A}(G_i, A)/D$. 
By the Scolem Theorem, $\Pi G_i/D_i$ is a faithful group over
the ultrapower $\Pi A/D$ of the ring $A$ over $D$. But there is 
the diagonal embedding of $A$ into $\Pi A/D$. So $\Pi G_i/D$ is 
a subgroup of the $A$--group $\Pi G_i^A /D$.\par
3) is a corollary of 1) and 2). \hfill $\Box$.

For quite a wide class of rings $A$ we can claim that torsion--free
nilpotent groups are faithful over $A$.

\begin{pr}
\label{5}
Let $G$ be a torsion--free nilpotent group, and $A$ an integral domain
of characteristic 0. Then $G$ is faithful over $A$.
\end{pr}

$\Box$ Let $A^\ast = A \otimes {\bf Q}$ be the tensor completion of
the ring $A$ over {\bf Q}. Then $A^\ast$ is a binomial domain
(i.e. it contains all binomial coefficients $(^a_n)$, $a \in A, \ 
n \in N$). It is known \cite{hall} that every torsion--free nilpotent group
$G$ has an $A^\ast$--completion in the sense of P. Hall. So $G$ 
is embeddible into an $A^\ast$--group, and hence into an $A$--group.
By remark~\ref{1}  $G$ is faithful over $A$.\hfill $\Box$

\begin{theorem}
Let $G$ be a residually torsion-free nilpotent group. Then $G$ 
is a faithful group over any integral domain $A$ of characteristic 0.
In particularly, free groups, free solvable groups and free polynilpotent 
groups  are faithful over any integral domain of characteristic 0.
\end{theorem}
$\Box$ The first statement follows from proposition 5  and the remark above.
Those concrete groups are residually torsion--free nilpotent: see, for example, \cite{mks} and \cite{hart}.
\hfill $\Box$ 

Finally, let us consider some examples of non--faithful groups.
\begin{example}
Let $G$ be a simple group containing an element of finite order. 
Then $G$ is not faithful over the field of rational numbers
${\bf Q}$. Moreover, $G^{\bf Q} = 1$. Indeed, any ${\bf Q}$--group has no
elements of finite order, so $G$ is not a subgroup of $G^{\bf Q}$. 
But in this case the homomorphism $\lambda: G \longrightarrow 
G^{\bf Q}$ has a nontrivial kernel. So the simplicity of $G$
implies that $\lambda(G)=1$, but $G^{\bf Q}$ is generated by $
\lambda(G)$. Therefore, $G^{\bf Q}=1$.
\end{example}

\begin{example}
Let $G$ be a torsion--free group with non--unique extraction of roots.
Then $G$ is not faithful over ${\bf Q}$, since $G^{\bf Q}$ has unique
extraction of roots. 
\end{example}

\section{Extensions of centralizers}
In the following we will use a construction of free extension of
centralizers of a  given group $G$. Let us describe this 
construction. \par

Let $G$ be a group. The centralizer of an element $v\in G$
is denoted by $C_G(v)$.
\par
If $G=\langle X \mid R \rangle$ is a presentation of $G$,
$Y$ is a set of words in $X$ and $t$ is a new letter
(not in $X$), then by $\langle G, \ t \mid [Y,t]=1 \rangle$
we will denote the group with the representation
$\langle X, \ t \mid R, [y,t] \rangle_{y \in Y}$.\par
Now we will define, perhaps, the simplest case of centralizer
extension.

\begin{opr}
The group $G(v,t)=\langle G,\ t \mid [C_G(v), t]=1 \rangle$
is called the direct of rank 1 extension of the centralizer
of the elememt $v$.
\end{opr}

It is easy to see that $G(v,t)$ can be obtained from $G$ by 
an HNN-extension with respect to the identity isomorphism
$C_G(v) \rightarrow C_G(v):$
$$
G(v,t) = \langle G,t \mid t^{-1}at=a, \ a \in C_G(v) \rangle,
$$
or as a free product with amalgamation:
$$
\langle G \ast (C_G(v) \times \langle t \rangle )\mid
C_G(v)=C_G(v)^\phi \rangle
$$
with amalgamation by the canonical monomorphism
$\phi : C_G(v) \rightarrow C_G(v) \times \langle t \rangle$.
\par

We will consider also the following much more general 
construction of  free centralizer extensions.

\begin{opr}
Let $G$ be a group, $C_G(v)=C$ the centralizer of an element 
$v$ from $G$, and $\phi :C\longrightarrow H$  a monomorphism
of groups such that $\phi (v)\in Z(H)$. Then the group
$$
G(v,H)=\langle G\ast H \mid C=C^\phi \rangle
$$
is called an extension of the centralizer $C_G(v)$ by the group
$H$ with respect to $\phi$.
\end{opr}
\par
Some particular cases of this construction are of special interest.
We will classify them with respect to the type of the embedding
$\phi: C \rightarrow H$. So the extension is 
 {\it central} iff $C^\phi \leq Z(H)$;
 {\it direct} iff $H=C^\phi \times B$;
{\it abelian} iff $H$ is an abelian group.

\begin{pr}
The following statements are true:
\begin{enumerate}
\item there exist the canonical embeddings $\lambda_G: G \longrightarrow
G(v,H)$ and \\ $\lambda_H: H \longrightarrow G(v,H)$;
\item $C_{G(v,H)}(v) \geq H$;
\item the group $G(v,H)$ has the following universal property:
for any group homomorphisms $\psi_G: G \longrightarrow N$ and 
$\psi_H: H \longrightarrow N$ compatible on $C$ and
$C^\phi$ there exists a unique homomorphism $\psi: G(v,H)
\longrightarrow N$ such that the following diagram is
commutative:

\medskip

\begin{center}

\begin{picture}(200,100)(0,0)
\put(0,100){$G$}
\put(85,100){$G(v,H)$}
\put(200,100){$H$}
\put(100,0){$N$}
\put(15,103){\vector(1,0){65}}
\put(195,103){\vector(-1,0){65}}
\put(15,93){\vector(1,-1){80}}
\put(105,93){\vector(0,-1){80}}
\put(195,93){\vector(-1,-1){80}}
\put(45,110){$\lambda_G$}
\put(155,110){$\lambda_H$}
\put(42,45){$\psi_G$}
\put(110,50){$\psi$}
\put(155,45){$\psi_H$}
\end{picture}

\end{center}

\medskip

\end{enumerate}

\end{pr}

\par
Now we will describe a construction which allows one to extend a
set of centralizers at once. 
\par
\begin{opr}
Let ${\cal C}=\{C_i = C_G(v_i) \mid i \in I \}$ be a set
of centralizers in the group $G$. Suppose $\phi_i: C_i 
\longrightarrow H_i$ is an embedding of $C_i$ into $H_i$
such that $\phi_i(v) \in Z(H_i), \ i \in I$. We can form a 
graph of groups

\medskip

\begin{center}

\begin{picture}(100,100)(0,0)
\put(0,50){$G$}
\put(100,100){$H_i$}
\put(100,50){$H_j$}

\put(103,75){\vdots}
\put(103,25){\vdots}
\put(15,57){\line(2,1){75}}
\put(15,53){\line(1,0){75}}
\put(50,85){$C_i$}
\put(50,60){$C_j$}
\end{picture}

\end{center}

\medskip

The fundamental group of this graph is called a tree extension of 
centralizers from ${\cal C}$ and is denoted by $G({\cal C, H}, \Phi)$,
where ${\cal H}=\{H_i \mid \ i \in I \},\ \Phi=\{ \phi_i \mid \ i \in I \}$.
\end{opr}
\par
Again, we will consider central, direct and
abelian extensions corresponding to the type of monomorphisms
from $\Phi$.
\par
By definition, the group $G({\cal C, \ H}, \Phi)$ is the union
(direct limit) of the chain of groups
$$
G=G_0 \hookrightarrow G_1 \hookrightarrow \cdots \hookrightarrow
G_\alpha \stackrel{\phi_\alpha}{\hookrightarrow} G_{\alpha+1}
\hookrightarrow \cdots
$$
where the set of centralizers ${\cal C}=\{C_\alpha \mid 
\alpha < \lambda\}$ is well ordered, $G=G_0,\  G_{\alpha+1}=
G_\alpha(v_{\alpha+1},H_{\alpha+1}),\ \phi_\alpha: G_\alpha
\longrightarrow G_{\alpha+1}$ is the canonical embedding from
proposition 5, and 
$G_\gamma=\lim_{\longrightarrow} G_\alpha$ (here $\alpha < \gamma$)
for a limit ordinal $\gamma$.\par
\begin{remark}
The properties of tree extensions of centralizers are similar to those of
ordinary extensions of centralizers: there exist the canonical embeddings
of $G$ and $H_i$ into $G({\cal C, H}, \Phi)$; the centralizer of $v_i$ in $G({\cal C, H}, \Phi)$ contains the group $H_i$ and the group $G({\cal C, H}, \Phi)$ has the corresponding universal property. Moreover, 
 groups $G({\cal C,H}, \Phi)$ for different well-orderings
have the same universal property, so they are isomorphic.   
\end{remark}

Finally, we introduce one more construction: the {\it iterated tree 
extension of centralizers.}
\par
Let $G$ be a group, ${\cal C}=\{C_G(v_i) \mid i \in I \}$ a set
of centralizers in $G$, and $\Phi = \{\phi_i: C_G(v_i) \rightarrow
H_i \}$ a set of embeddings of these centralizers. Then one can
construct the group $G_1=G({\cal C, H}, \Phi)$ and consider again
the set of centralizers ${\cal C}_1$ in $G_1$, the set of embeddings
$\Phi_1$, and the corresponding set of groups ${\cal H}_1$. As above,
the group $G_2 = G_1({\cal C}_1, {\cal H}_1, \Phi_1)$ can be
constructed. We can repeat this process using induction on ordinals:
if a group $G_\alpha$ has been constructed for some ordinal 
$\alpha$ and the set of centralizers ${\cal C}_\alpha$ in $G_\alpha$,
the set of groups ${\cal H}_\alpha$ and the corresponding set of
embeddings $\Phi_\alpha$ have been chosen, then 
$$
G_{\alpha+1} = 
G_\alpha({\cal C}_\alpha, {\cal H}_\alpha, \Phi_\alpha),
$$ 
$$
G_\gamma = \lim_{\stackrel{\longrightarrow}{\alpha <\gamma}} G_\alpha,
$$ 
for a limit ordinal $\gamma$.
Thus,
for any ordinal $\delta$, starting from a given group $G$, we can obtain
the group $G_\delta$ by iterating the process above.

\begin{opr}
The group $G_\delta$ is called an iterated tree extension of centralizers
up to level $\delta$.
\end{opr} 

\begin{remark}
Any iterated tree extension of centralizers can be obtained by
transfinite induction using operations of extension of centralizers
and direct limits. The group $G$ is embeddible in $G_\delta$,
and $G_\delta$ has a universal property similar to that in proposition 5
\end{remark}

\section{Canonical forms of elements in groups which are obtained
by abelian extension of centralizers}

In this section we will consider canonical and reduced forms of
elements of groups constructed by extensions of centralizers.
Of course, there are general theorems about canonical forms
in amalgamated free products, but we need some special concrete forms
corresponding to this particular situation.\par

Let $G$ be a group, $v \in G$ an element such that $C_G(v)=C$ is 
a maximal abelian subroup of $G$,
and $\phi: C \hookrightarrow A$ an embedding of $C$ into some
fixed abelian group $A$. In exponential groups it will be
convenient to use an exponential writing for elements of $C$
and $A$. Namely, $A \simeq v^A = \{v^a \mid a \in A \}$,
$C \simeq v^C = \{v^c \mid c \in C \}$, in particular,
$\langle v \rangle \simeq v^{\bf Z} = \{v^n \mid n \in {\bf Z} \}$.
Consider the group $G^\ast$ which is obtained from $G$ by 
$A$-extension of the centralizer $C$:
$$
G^\ast=G(v,A)=\langle G\ast v^A \mid v^C=v^{C^\phi} \rangle
$$
with respect to the embedding $\phi$.
\par \medskip
Denote by $S$ some system of right representatives of $A$ by
$C$. Then $v^S$ is a system of right representatives of 
$v^A$ by $v^C$. Fix any system ${\cal R}$ of right representatives
of $G$
by $v^C \simeq C.$ For $g \in G$ we denote by
$\bar{g}$ the representative of $g$ in this fixed system of 
representatives. Clearly, $g=v^c \bar{g}$ for some $c \in C$.
From general theory of free products with amalgamation 
we know that every element $g \in G^\ast$ can be written in the
form
\begin{equation}
g=g_1v^{a_1}\cdots g_nv^{a_n}g_{n+1} \label{eq1}
\end{equation}
$a_i \in A, \ g_i \in G$. Moreover,
any element $g \in G^\ast$ admits a unique 
decomposition of the following type, which is called {\it a canonical form
of g}:
$$
g=hg_1v^{s_1}\ldots g_nv^{s_n}g_{n+1} 
$$
where $h \in v^C, \ g_i \in {\cal R}, 
\ s_i \in S,\ s_i \neq 0$ and $g_i \neq 1$ (except, maybe,
$g_1,\ g_{n+1}$). We suppose here that $g$ itself is a canonical form
for $g \in G$.
\par \medskip
Any product of type (\ref{eq1}) (or, more precisely, any sequence
$(g_1, v_1^{a_1}, g_2, \ldots, v^{a_n}, g_{n+1})$)
can be reduced to a canonical
form by the following {\it rewriting process}:
\begin{itemize}
\item
if $a_n \in C$, then
$$
g_1v_1^{a_1} \cdots g_nv^{a_n}g_{n+1}
\rightarrow g_1v^{a_1} \cdots g_{n-1}v^{a_{n-1}}(g_nv^{a_n}g_{n+1}),
$$
so we get a shorter decomposition of $g$ and can apply induction;
\item if $a_n \not \in C$, then 
$$
g_1v^{a_1} \cdots g_nv^{a_n}g_{n+1}
\longrightarrow g_1v^{a_1} \cdots (g_nv^{k+t})v^{s_n} \overline{g}_{n+1},
$$
where $a_n=k+s_n, \ g_{n+1} = v^t \overline{g}_{n+1}, \ k,t \in C,\ 
s_n \in S$;

\item if $(g_nv^{k+t}) \in C$,
then 
$$
g_1v^{a_1}\ldots v^{a_{n-1}}(g_nv^{k+t})v^{s_n}\overline{g}_{n+1}
\longrightarrow g_1v^{a_1}\ldots 
(g_{n-1}g_nv^{k+t})v^{a_{n-1}+s_n}\overline{g}_{n+1},
$$
and we get a shorter decomposition of $g$. After that the process is
performed by induction.
\item if $g_nv^{k+t} \not \in C$, then one can apply the 
rewriting process to $g_1v^{a_1} \cdots v^{a_{n-1}}g_nv^{k+t}$,
keeping the tail unchanged.
\end{itemize}
\par
It is easy to see that these transformations
carry $g$ to the canonical form.

\par
It is more convenient to use the reduced form of elements of $G^\ast$.
By definition, an element $g$ is in {\it reduced form }
\begin{equation}
g\ =\ g_1v^{a_1}\ldots g_nv^{a_n}g_{n+1} \label{eq2}
\end{equation}
if $a_i \not\in C, i=1,\ldots ,n, g_j\not\in C, j=2,\ldots ,n$.
\par

From the rewriting process it is clear that a
 form of type (1) is reduced iff the number $n$ is
minimal through all decompositions of $g$ of type (1).

\begin{pr}
Let
$$
g=u_1v^{a_1}\ldots u_mv^{a_m}u_{m+1}
$$
$$
g=w_1v^{b_1}\ldots w_nv^{b_n}w_{n+1}
$$
be two reduced forms of an element $g$. Then: 
\begin{enumerate}
\item $m=n$ and $a_1=b_1+k_1,\ldots ,a_m=b_m+k_m$ for some $k_i \in C$;
\item $u_{m+1}  =  v^{t_m}w_{m+1},\ u_m=v^{t_{m-1}}w_mv^{-t_m-k_m},\ \
\ldots, u_1= w_1v^{-t_1-k_1}$ for some $t_j\in C, j=1,\ldots ,m$.
\end{enumerate}
\end{pr}
$\Box$ By the rewriting process. $\Box$
\par
For a direct extension of a centralizer $C$ by a group $H$
(i.e. when $A = C \oplus H$), the group $H$ is itself a system
of representatives of $A$ by $C$ ($H=S$).
 In this case it is 
very convenient
to use semicanonical forms of elements.\par 
By definition, a reduced
form (\ref{eq2})
$$
g=g_1v^{s_1} \cdots v^{s_n}g_{n+1}
$$
is {\it semicanonical}, iff $1\not=s_i \in S$ for all $i$.\par
From proposition 6 one can easily deduce the following
\begin{co}
Any two semicanonical forms of $g$ can be transformed from one to another
by a finite sequence of commutations of the form $v^sv^t=v^tv^s$,
where $t\in C,\ s\in S$.
\end{co}

Suppose now that the notions of canonical, semicanonical
and reduced forms of elements of the group $G$ are already 
introduced. In this case one can extend them to the group $G^\ast$.\par
We will say that canonical (semicanonical, reduced) forms of an element
$$
g=hg_1v^{a_1} \cdots v^{a_n}g_{n+1}
$$
{\it agree with those on $G$}, iff the elements
$h, g_1, \ldots, g_{n+1}, v$ are in canonical (semicanonical,
reduced) forms in the group $G$. \par
Now we are able to introduce the corresponding forms of elements
of tree and iterated tree extensions of centralizers.\par
Let $G^\ast$ be the union of a chain of groups
$$
G = G_0 \hookrightarrow G_1 \hookrightarrow \cdots 
\hookrightarrow G_\alpha \hookrightarrow \cdots
$$
where every $G_{\alpha + 1}$ is obtained from $G_\alpha$
by extension of centralizers, and $G_\gamma = \lim_{\rightarrow}
G_\alpha$ for a limit ordinal $\gamma$. Then the notions of
canonical (semicanonical, reduced) forms are introduced on $G_{\alpha +1}$
agreeably with those on $G_\alpha$. For a limit ordinal $\gamma$ 
these forms are inherited form all the previous terms of the chain.
Then {\it canonical (semicanonical, reduced) forms} on the group
$G^\ast$ are well--defined.

\section{Groups with conjugately separated maximal abelian subgroups}
\par
In this section we introduce the class of CSA$^\ast$--groups 
which will play a 
fundamental role in our investigation of $A$-free groups, tensor
completions and groups with length functions. This class is quite
wide; for example, it contains all torsion-free hyperbolic groups, 
all groups acting freely on $\Lambda$--trees, and all universally
free groups.
\begin{opr}
A subgroup $H$ of a group $G$ is called {\it conjugately separated} 
(or {\it   malnormal}) if 
$H \cap H^x = 1$ for any $x \in G \setminus H$.
\end{opr}

\begin{opr}
A group G is called a CSA-group if all its maximal abelian subgroups
are conjugately separated.
\end{opr}

\par
It is clear that every abelian group is a CSA--group.\par
The following proposition gives us some intuition about CSA-groups.

\begin{pr}
\label{3.1}
Let $G$ be a CSA-group. Then the following statements are true:
\begin{enumerate}
\item (SA--property) if $M_1$ and $M_2$ are distinct maximal abelian 
subgroups of $G$, then $M_1 \cap M_2 = 1$;
\item any maximal abelian subgroup of $G$ coincides with its
normalizer. 
\end{enumerate}
\end{pr}
\par
$\Box$ 

1) If $1 \neq x \in M_1 \cap M_2$ and $y \in M_2 \setminus M_1$, then
$x \in M_1^y \cap M_1$, so $M_1^y \cap M_1 \neq 1$.
\par \medskip
2) is evident.
$\Box$

\begin{opr}
A group $G$ is called 
commutative transitive if the relation ``$a$ commutes with $b$'' is
transitive on the set $G \setminus{1}$;
\end{opr}

\begin{remark}
The following conditions are equivalent to one another:
\begin{enumerate}
\item a group $G$ has SA-property;
\item a group $G$ is commutative transitive;
\item in a group $G$ the centralizer of any nontrivial
element is abelian.
\end{enumerate}
\end{remark}

\begin{remark}

Proposition 7 shows that the class CSA is a 
subclass of the class of groups with transitive commutation. The following example
proves that CSA is a proper subclass of these classes.
\end{remark}
$\Box$
Let $G$ be a free metabelian group. Then all proper centralizers of $G$ are
commutative. But $G$ has a
normal abelian subgroup, so it is not a CSA--group. \hfill $\Box$
\par \medskip
In fact, the transitivity of commutation is a very strong property 
and implies some other interesting features. We collect them in the
following
\begin{pr}
\label{3.2}
Let $G$ be a group such that commutation is an equivalence relation
on $G \setminus \{1\}$. Then
\begin{enumerate}
\item the centralizer of any nontrivial element of $G$ is a maximal 
abelian subgroup of $G$; conversely, any maximal abelian subgroup of 
$G \neq 1$ is a centralizer of some nontrivial element;
\item if $G$ is torsion--free, then for $x,y \in G$
 $x^n=y^m \Longrightarrow [x,y]=1;$\\ 
 $x^n=y^n \Longrightarrow x=y;$ 
\item if $G$ is nonabelian, then its center is trivial;
\item if $G$ is nonabelian, then it is indecomposable into a 
nontrivial direct product.
\end{enumerate}
\end{pr}
$\Box$ 1), 3) and 4) are quite evident. The second implication 
in 2) is a corollary of the first one. Let us prove the first claim
of 2).
\par 
Suppose $x^n=y^m$. Then the neighboring terms in the sequence 
$x, x^m, y^n, y$ are nontrivial and commute pairwise; by transitivity
of commutation we obtain $[x,y]=1.$ \hfill $\Box$

\begin{remark}
$G$ is a CSA-group iff all its  centralizers of nontrivial
elements are abelian
and conjugately separated (see statement 1 of 
proposition 8).
\end{remark}

Thus, CSA is a property of centralizers. Let us consider a typical
situation of conjugate unseparability. A cyclic subgroup 
$\langle y \rangle = M < G$ is not conjugately separated iff
there is some $x \in G$ such that $M^x \cap M \neq 1$, or,
equivalently, $x$ and $y$ satisfy the identity
$x^{-1}y^rx=y^s$ for some integers $r$ and $s$. One-relator
groups with this relation have a special name 
and play a crucial role in what follows. Groups of the form
$$
G_{r,s} = \langle a,b \mid a^{-1}b^ra=b^s \rangle, \ r,s \in {\bf Z}, 
r > 0,
$$
are called {\it Baumslag -- Solitar groups} \cite{baso}. We will need the 
following properties of these groups:
\begin{itemize}
\item if $r=s=1$, then $G_{r,s} \simeq {\bf Z} \times {\bf Z}$ is
an abelian group;
\item if $r=1$ and $s \neq 1$, then $G_{r,s}$ is nonabelian solvable
and for $s \neq -1$
$$
G_{1,s} = \langle a,b \mid a^{-1}ba = b^s \rangle \simeq 
{\bf Z}[\frac{1}{s}] >\hspace{-.08in}\lhd {\bf Z}
$$
where ${\bf Z} \simeq \langle a \rangle,  $
${\bf Z}[\frac{1}{s}] \simeq {\rm n.cl}(b) $ - the normal closure of 
$\langle b \rangle$; \par
for $s= -1$ 
$$
G_{1,s} = \langle a,b \mid a^{-1}ba = b^{-1} \rangle \simeq
\langle b \rangle >\hspace{-.08in}\lhd \langle a \rangle.
$$
\item if $r > 1, \ \mid s \mid >1$, then $G_{r,s}$ is a nonsolvable
group and $[a^{-1}ba,b] \neq 1$ (see, for example, \cite[p.281]{mks}).
The case $r >1, s = \pm 1$ coincides with the second one after renaming
the generators.
\end{itemize}
\par
In the following proposition we collect some properties of 
CSA--groups.
\begin{pr}
Let $G$ be a nonabelian CSA--group. Then
\begin{enumerate}
\item  centralizers of nontrivial elements of $G$ are abelian
and coincide with maximal abelian subgroups of $G$ (in particular,
commutation is transitive);
\item $G$ has no nontrivial normal abelian subgroups;
\item $G$ has a trivial center;
\item $G$ is directly indecomposable;
\item $G$ has no non--abelian solvable subgroups;
\item $G$ has no nonabelian Baumslag -- Solitar subgroups;
\item if $G$ is torsion--free, then $G$ is a D--group, i.e.
the operations of extracting roots in $G$ are defined 
uniquely.
\end{enumerate}
\end{pr}
$\Box$ 1) is follows from item 3) of proposition 7 and item 1)
of proposition 8. 2) is clear, and 3) is a corollary of 2).
4) follows from 1) and item 4) of proposition 8. In proposition 10
 below it
will be proved that a subgroup of a CSA--group is again a CSA--group. 
So to prove 5) one only needs to prove that a non--abelian solvable
group is not a CSA--group. But the last statement follows from 2). To prove 6),
it suffices to prove that a non--abelian Baumslag--Solitar group $G_{r,s}$
is not a CSA--group. If $r=1$ and $s \neq 1$, then $G_{r,s}$ is
nonabelian solvable (see remarks before the proposition), and hence by
5) it is not a CSA--group. If $r > 1$ and $ \mid s \mid >1$, then
$[a^{-1}ba,b] \neq 1$. Suppose now that $G_{r,s}$ is a CSA--group, then
commutation on nontrivial elements of $G$ must be transitive. In the
sequence $a^{-1}ba, \ a^{-1}b^ra, \ b^s\ , b$ all neighbours commute,
but the first and the last elements do not - a contradiction with
the transitivity of commutation. 7) follows from item 3) of proposition 7 
and item 2) of proposition 8. \hfill $\Box$ 

\begin{pr} 
\par
\begin{enumerate}
\item The class of groups CSA is finitely universally axiomatizible.
\item The class CSA is closed under ultraproducts.
\item The class CSA is closed under taking subgroups.
\end{enumerate}
\end{pr}

\par \medskip
$\Box$  To prove 1) it suffices to write down the conditions of
commutativity of centralizers and their conjugate separation by 
universal formulas. This is straightforward.
2) follows from the fact that this class is axiomatizible.
3) follows from the universal axiomatization of this class.
\hfill $\Box$
\par \medskip
Free groups are principal examples of nonabelian CSA-groups. 
The following results show that the class CSA contains almost all
groups which are ``close'' to free groups.

\begin{pr} 
\label{3.6}
Let $G$ be a torsion-free group such that
 the centralizer of any nontrivial element is cyclic.
Then $G$ is a CSA-group.
\end{pr}

\par \medskip
$\Box$ Let $G$ satisfy conditions of Proposition 11. Then proper centralizers of
$G$ are abelian (in particular, commutation is transitive on 
$G \setminus \{1\}$), and by remark 4 we only need to prove that these
centralizers are conjugately separated. Consider any maximal abelian subgroup 
$M$ of $G$; $M=\langle y \rangle$. Suppose that for some $x \in 
G \setminus M$ we have $M \cap M^x \neq 1$. This means that
$x^{-1}y^sx=y^r$ for some integers $s,r$. Let us denote $z=x^{-1}yx$.
Then $z^s = y^r$ and by proposition 8 $[z,y]=1$. So $z \in 
\langle y \rangle$ and $x^{-1}yx = y^p$ for some integer $p$. 
Therefore the subgroup $\langle x,y \rangle$ is a quotient group
of the Baumslag -- Solitar group $G_{1,p} \simeq \langle a,b \mid
a^{-1}ba=b^p \rangle$ (see remarks before proposition 9)
under the epimorphism $\phi: a \mapsto x,\  b \mapsto y$. We
know that $G_{1,p} = {\rm n.cl} \langle b \rangle >\hspace{-.08in}\lhd 
\langle a \rangle$, where n.cl$\langle b \rangle \simeq {\bf Z}[\frac{1}{p}]$,
$\langle a \rangle \simeq {\bf Z}$. Suppose that $ b_1a^n \in 
{\rm ker}\phi$, i.e. $y_1x^n=1, \ y_1 \in {\rm n.cl}(y)$. Then
in the sequence $y,y_1, x^n, x$ neighbours commute, and hence by
transitivity of commutation $[y,x]=1$ (but this contradicts
to our choice of $x$ and $y$) or else $y_1$ or $x^n$ is trivial.
The last case implies that $y_1=1, \ n=0$. Hence $\langle x,y
\rangle \simeq G_{1,p}$ - a contradiction, because
${\bf Z}[\frac{1}{p}]$ is a subgroup of $G_{1,p}$ and is not cyclic.
\hfill $\Box$

\begin{pr} 
\label{3.7}
\par\medskip
\begin{enumerate}
\item Any torsion-free hyperbolic group is a CSA-group.
\item Any group acting freely on a $\Lambda$-tree is a CSA-group.
\item Any universally free group is a CSA-group.
\end{enumerate}
\end{pr}

$\Box$ It is known (see, for example, \cite{gr}) that the centralizer 
of any element
of infinite order in a hyperbolic group $G$ has a cyclic subgroup of finite
index (i.e. the stabilizer is virtually cyclic). So we only need to prove
that any torsion-free virtually cyclic group $H$ is infinite cyclic. It is
clear that $H$ is finitely generated and all subgroups of $H$ are also
virtually cyclic. By induction on the number of generators of $H$, we may
assume that $H$ is 2-generated, $H= \langle a,b\rangle$ and moreover 
$\langle b\rangle\lhd H,\ a^m\in \langle b\rangle,\ a^{-1}ba=b^n$,
for some $n,m$. This implies $b=a^{-m}ba^m = b^{n^m}$. Therefore 
$1=n^m$.
Hence $n=\pm 1$ and $[a,b]=1$
or $[a^2,b]=1$. In the first case $H$ is abelian torsion--free and
virtually cyclic, therefore $H$ is a cyclic group. In the second case 
$a^2$ is a central element of $H$, and so is the element $a^{2m}=b^k$
for some integer $k$. But $a^{-1}b^ka=b^{-k}$, and hence $b^k=b^{-k}$,
which contradicts to the torsion--free condition on $G$.\par

Statement 2) is proved in \cite[corollary 1.9]{bass}. \hfill $\Box$
\par
\medskip
In \cite{vnr} it was proved that every universally free group
acting freely on some $\Lambda$-tree.
Now 3) follows from 2).\hfill $\Box$

\section{Preservation of the CSA--property under free products,
direct limits, and abelian extensions of centralizers}

In this section we will prove that the CSA--property is preserved
under free products and abelian extensions (tree extensions or
iterated tree extensions) of centralizers for groups without elements
of order 2 (non-abelian CSA groups automatically do not contain elements of order 2). In any case it is preserved under direct limits.

\begin{opr}
The subclass CSA$^\ast$ of CSA consists of all groups without elements
of order 2.
\end{opr}

\begin{remark}
Any non-abelian CSA-group is also CSA$^\ast$-group.
\end{remark}
$\Box$ If $G$ is non-abelian CSA-group and $u^2=1$ for some nontrivial 
$u \in G$, then the normal closure of $u$ in $G$ is also non-abelian 
(otherwise it would contradict to CSA-property). Hence, there are two 
noncommutative involutions in $G$, they generate non-abelian metabelian 
subgroup of $G$ - contradiction with CSA-property.
$\ \ \ \ \ \Box$
\par \medskip

\begin{theorem}
The class CSA$^\ast$ is closed under free products.
\end{theorem}
$\Box$ Let $G_1$ and $G_2$ be CSA$^\ast$--groups and  $G=G_1\ast G_2$.
First of all, $G$ has no elements of order two if $G_1$ and $G_2$
do not (see \cite[cor. 4.4.5]{mks}).  By the Kurosh Subgroup Theorem
a maximal abelian subgroup $M$ of the group $G$ has one of
the following two types:
\begin{enumerate}
\item[a)] $M \leq G_i^g$ for some $g \in G, \ i=1,2$.
\item[b)] $M= \langle z \rangle$ is an infinite cyclic group generated 
by an element $z$ of cyclically reduced length $\| z \| \geq 2$. 
\end{enumerate}
By definition,
$\| z \| = {\rm min}\{ \mid y\mid \ \mid y$ is cyclically reduced and 
conjugate to $z \}$. Clearly, $\|z\| \geq 2$ iff $z \not \in G_i^g$
for any $g \in G, \ i=1,2$.
We claim that all centralizers of nontrivial elements in $G$ are abelian and, consequently,
commutation in $G \setminus \{1\}$ is transitive. Indeed,
let  $1 \neq a \in G, \ b, c \in C_G(a)$. There are maximal
abelian subgroups $M, N < G$ with $a,b \in M$ and $a,c \in N$.
We have the following cases:
\begin{enumerate}
\item $M < G_i^x, N < G_i^y$, $i=1,2, \ x, y \in G$. Then
$G_i \cap G_i^{xy^{-1}} \neq 1 $, hence $xy^{-1} \in G_i$, and
$G_i^x = G_i^y$. So $a,b,c \in G_i^x$. But since $G_i$ is a 
CSA--group, it follows that $[b,c]=1$.
\item The case $M < G_i^x, N < G_j^y, i \neq j$ is impossible. 
\item The case $M= \langle z \rangle$ with $ \| z \| \geq 2, \ N < G_i^x$
is impossible, because $a \in M$ implies $\| a \| \geq 2$, and 
on the other hand, $a \in N$ implies $\| a \| = 1$.  
\item $M = \langle z \rangle, N = \langle y \rangle, \| z \| \geq 2,
\| y \| \geq 2$. Then $a=z^n=y^l$. So $y^l \in Z(H)$, where
 $H=\langle z, y \rangle$. By the Kurosh subgroup theorem, $H$ is a free
product $H=F \ast (\ast G_i^\alpha)$, where $F$ is a free group.
Since $H$ has a nontrivial centre, it equals to one of the factors. 
If $H=F$, then $H$ is cyclic, and $[z,y]=1$. Therefore $[b,c]=1$. 
The case $H < G_i^\alpha$ is impossible because $\| z \| \geq 2$.
\end{enumerate}
This proves that $C_G(a)$ is abelian. 
\par
 Suppose that $M \leq G_i^g$, then $M$
is conjugately separated iff its automorphic image $M^{g^{-1}}$ is
conjugately separated. Conjugating by $g^{-1}$ if necessary, we
may assume that $M \leq G_i$. It is a well-known property of free 
products (see \cite[cor. 4.1.5]{mks}) that if $x$ is not in $G_i$,
then $G_i^x \cap G_i = 1$. Hence if $M^x \cap M \neq 1$, then 
$x \in G_i$, but $G_i$ is a CSA-group, therefore $x \in M$.
Now suppose
 that $M=\langle z \rangle$ and $\| z \| \geq 2.$ If
$M^x \cap M \neq 1$, then $x^{-1}z^nx=z^m$ for some $n$ and $m$.
Hence $\| z^n\| = \| z^m \|$. But it is easy to see that
$\| z^n \| = \mid n\mid \cdot \| z \|$, and the equality above implies
that $\mid m\mid =\mid n\mid$. In the case $n=m$ we have $[x, z^n]=1$, and
by transitivity of commutation  $[x,z]=1$. Consequently, $x \in 
\langle z \rangle$.  Assume now that $n=-m$. In this case 
$x^{-1}z^nx=z^{-n}$, which implies $x^{-2}z^{n}x^2=z^n$. 
So $[x^2, z^n]=1$, and by transitivity $[x,z]=1$, as above.
$\ \ \ \ \ \Box$
\par \medskip
The following theorem plays the main part in the rest of the paper.

\begin{theorem}
The class CSA$^\ast$ is closed under extensions of centralizers
by abelian groups without elements of order 2.
\end{theorem}
$\Box$ Let $G$ be a CSA$^\ast$--group, $C=C_G(v)$ the centralizer of an
element $v \in G$, and $\phi: C\hookrightarrow A$ an embedding of $C$ into
an abelian group $A$ which has no elements of order two. We need to prove
that the extension $G(v,A)$ is also a 
CSA$^\ast$--group. The group $G(v,A)$ has no elements of order
two (\cite[cor. 4.4.5]{mks}). Therefore it remains to check the CSA--property.
Before proving this, we establish a lemma
which describes centralizers in the group $G(v,A)$.

\begin{lemma}
Let $x \in G(v,A)$. Then:
\begin{enumerate}
\item if $x \in G^g$, then $C(x) \leq G^g$;
\item if $x \in A^g$, then $C(x) = A^g$;
\item if $\| x \| \geq 2$, then $C(x)=\langle z \rangle $ and $\| z \| \geq 2$.
\end{enumerate}
\end{lemma}

$\Box $ We will use the following description of commuting elements in free products
with amalgamation (see \cite{mks}): $[x,y]=1$ iff one of the following 
conditions holds:
\begin{itemize}
\item $x$ or $y$ belongs to some conjugate of the amalgamated subgroup
$C=C_G(v)$;
\item neither $x$ nor $y$ is in a conjugate of $C$, but $x$ is in a
conjugate of a factor ($G$ or $A$), then $y$ is in the same conjugate of 
the factor;
\item neither $x$ nor $y$ is in a conjugate of a factor, then 
$x=g^{-1}cgz^{n}, \ y=g^{-1}c^\ast gz^m$, where $c, c^\ast \in C$,
and $g^{-1}cg, \ g^{-1}c^\ast g$ and $z$ commute pairwise.
\end{itemize}

Let us consider these three possible cases.
\par \medskip
Suppose that $x \in C^g, \ y \in G(v,A),$ and $[x,y]=1$. We can assume
(conjugating by $g^{-1}$ if necessary) that $x \in C$. Let 
$y=cg_1a_1 \ldots g_na_n$ be the canonical form of $y$, where
$c \in C$, $g_i$ are right representatives of $G$ mod $C$, 
 $a_i$ are right representatives of $A$ mod $C$, and $g_n \neq 1$.
Then $xy=yx$ implies
$$
xcg_1a_1 \ldots g_na_n = cg_1a_1 \ldots g_na_nx = cg_1a_1 \ldots g_nxa_n
$$
Hence $g_nx=x_1g_n$ for some $x_1 \in C$ or $g_nxg_n^{-1} = x_1$.
But $G$ is a CSA-group and $C$ is its nontrivial centralizer  --
a contradiction. Therefore if $[x,y]=1$, then $y=ca_n \in A$.
And we have case 2) of the lemma. 
\par
Suppose that $x$ belongs to a conjugate of a factor $G$ or $A$, but not to 
$C$. Then $y$ is in the same factor, and we have case 1) or 2) of the lemma. 
\par
Finally, suppose that $\| x \| \geq 2$ and $[x,y]=1$ for $y \in G(v,A)$.
Then $x=g^{-1}cgz^n, \ y=g^{-1}c^\ast g z^m, \ \| z \| \geq 2$, and
$z, g^{-1}cg$ and $g^{-1}c^\ast g $ commute in pairs. Conjugating by 
$g^{-1}$, we reduce the problem to the case  $x=cz_1^n, y=c^\ast z_1^m,
z_1=gzg^{-1}$. If $c \neq 1$ or $c^\ast \neq 1$, then $z_1 \in A$
according to the proof above, but this contradicts the condition
$\| z \| \geq 2$. So $c=1, \ c^\ast=1$ and we have the case 3) 
of the lemma. $\ \ \ \ \ \Box$

\begin{co} 
\begin{enumerate}
\item All centralizers of nontrivial elements of the group 
$G(v,A)$ are abelian, and the group $G(v, A)$ is commutative transitive;
\item The cenralizer of the element $v$ in the group $G(v,A)$ is exactly
the subgroup $A$.
\end{enumerate}
\end{co}

\begin{co}
A maximal abelian subgroup of $G(v,A)$ can be of one of the following
three types:
\begin{enumerate}
\item $M^g$, where $M$ is a maximal abelian subgroup of 
$G$, $M \neq C^f$ ($M$ is not conjugate to $C$) and $g \in G(v,A)$;
\item $A^g, \ g \in G(v,A)$;
\item $\langle z\rangle,$ where $\| z\|\geq 2.$
\end{enumerate}
\end{co}

The following lemma describes  conjugated
elements in the group $G(v,A)$.
\begin{lemma} Let $g$ be a cyclically reduced element of $G(v,A)$.
Then:
\begin{enumerate}
\item[a)] if $g$ is conjugate to an element $h$ in some factor, 
then $g$ and $h$ are in the same factor and are
conjugate in it;
\item[b)] if $g$ is conjugate to a cyclically reduced element 
$p_1\cdots p_n$, where $n\geq 2$, then $g$ can be obtained by cyclically
permuting $p_1,\ldots p_n$ and then conjugating by an element from $C$.
\end{enumerate}
\end{lemma}

$\Box$ There is a description of conjugated elements in arbitrary free
products with amalgamation \cite{mks}. In a general situation, one has our cases a),
b) and one more. Namely, $g$ is conjugate to some $c \in C$. In this case, 
according to the description in \cite{mks}, there are elements $c_i \in C$ such that 
in the sequence $c, c_1, \ldots, c_n, g$ consecutive terms are conjugate
in a factor. In our case the factor $A$ is abelian, therefore $g$ and $c$
are conjugate in $G$. \hfill $\Box$ \par
{\it Proof of the theorem.} 
Let $N=M^g$ be a maximal abelian subgroup of $G$ of the type 1). We
need to prove that $N^x \cap N \neq 1 \Longrightarrow \ x \in  N$
for any $x \in G(v,A)$. Conjugating if necessary, we can assume
that $N=M \leq G$. Suppose that there exist nontrivial elements $f,h \in M$
such that $f^x=h$.  
According to statement a) above, $f^y=h$ for some $y \in G$. 
But then $M \cap M^y \neq 1$, and consequently $y \in M$ and $f=h$
(because $G$ is a CSA--group). 
This implies that $[f,x]=1$ and $x \in G$. Hence $x \in M$. \par
Now let $N=A^g$, $g \in G(v,A)$. Again, we can consider only the case $N=A$.
Suppose that $A \cap A^x \neq 1$ for some $x \in G(v,A)$. Then 
$f^x=h$ for some $f,h \in A$. If $f$ or $h$ is not in $C=C_G(v)$,
then by lemma 2 a) $f$ and $h$ are conjugate in $A$, but $A$ is abelian,
so $f=h$. If both $g$ and $h$ are in $C$, then they are both in $A$ and $G$,
and we can apply a) as before.
\par
Let $N=\langle z\rangle$ be an infinite cyclic group, $\|z \| \geq 2$.
Suppose that $N\cap N^x\neq 1$. Then $(z^n)^x=z^m$ for some integers 
$m$ and $n$. Now we are in the case b) of the description.
To finish the proof in this case, one can repeat word by word the proof
from the last part of theorem 4.
\hfill $\Box$ 

\begin{theorem}
The class CSA is closed under direct limits.
\end{theorem}
$\Box$ Let ${\cal G}=\{G_i, \phi_j^i, I\}$ be an inductive
system of CSA--groups $G_i$, $i \in I$, and $G^\ast=\lim_{\longrightarrow} 
G_i$.
The following lemma describes centralizers of nontrivial elements in $G^\ast$.
\begin{lemma}
Let $G^\ast=\lim_{\longrightarrow} G_i$ be a direct limit of
CSA--groups $G_i$, $i \in I$.
If $M$ is a maximal abelian subgroup of $G^\ast$, then for
any $i \in I\ M_i =G_i\cap M$ is either maximal abelian 
in $G_i$ or trivial. Moreover, 
$ M = \lim_{\longrightarrow} M_i$
\end{lemma}
$\Box$ 
First of all, commutation is transitive on $G^\ast \setminus \{1\}$.
Indeed, any three elements $x,y,z\in G^\ast$ are in some $G_i$,
but $G_i$ is a CSA--group. \par
Let $M_i$ be nontrivial. If $M_i$
is not maximal in $G_i$, then $[M_i ,x]=1$ for some $x\not\in
M_i$. Let us prove that $x\in M$. By transitivity of
commutation, $M$ is a centralizer of
some element $z\in G$. There exists $k$ such that $i < k$
and $x,z\in G_k$. From transitivity of commutation 
we deduce that $[x,z]=1$, and therefore $x\in M$, which contradicts
to the condition $x\not\in M_i$. The equality
$M = \lim_{\longrightarrow} M_i$
is evident. \hfill $\Box$
\par \medskip
Let us continue the proof of theorem 6.\par
 Suppose that $M$ is a maximal
abelian subgroup of $G^\ast$, and for some $f\in G^\ast\ M\cap
M^f\neq 1$. Then $x=y^f$ for some $x,y\in M$. We can choose 
$k \in I$ such that $x,y,f\in G_k$ and therefore
$M_k\cap M_k^f\neq 1$ in $G_k$. Hence $f\in M_k
<M$. This proves that $G^\ast$ is a CSA--group. 
\hfill $\Box$
\begin{theorem}
The class CSA$^\ast$ is closed under tree abelian extensions of
cenralizers and iterated tree abelian extensions of centralizers
by abelian groups without elements of order $2.$
\end{theorem}
$\Box$ By definition, a tree (iterated tree) extensions of centralizers
of a group $G$ is constructed by transfinite induction on ordinals:
$G=G_0<G_1<\cdots <G_\alpha <\cdots.$ Following this construction,
it suffices to prove by induction that every term of this chain
lies in the class CSA$^\ast$ if the initial group $G$ belongs to
CSA$^\ast$. For a nonlimit ordinal, this was proved in theorem 5,
for a limit ordinal in theorem 6. \hfill $\Box$

\section{Tensor completion of CSA--groups}
In \cite{mr1}, using general category theory, we proved that for any group
$G$ and any ring $A$ the tensor completion $G^A$ exists, but 
we did not point out any idea of what this construction is explicitly.
\par  
In this section we, following the method of G. Baumslag \cite{bau}, 
describe the explicit construction of a tensor completion
of an arbitrary CSA$^\ast$--group $G$ using extensions of centralizers.
In fact, $G^A$ appears as an iterated tree extension of centralizers in $G$
for a very wide class of rings $A$. In particular, $G^A$ is again a 
CSA$^\ast$--group, and we obtain a lot of information about it.
Moreover, we can apply the technique of amalgameted free products to
investigate $G^A$. Thus, one has the canonical and reduced forms of
elements in $G^A$, descriptions of commuting and conjugating elements, 
elements of finite order and so on.
\par
\medskip
We will construct the $A$--completion of a group $G$ step by step,
on each stage obtaining partial $A$--groups.
 \par \medskip

{\it Abelian case}
\par\medskip
 In the
abelian case the construction of tensor completion is very similar
to ordinary tensor multiplication by the ring $A$ (one only
needs to make the existing partial action of $A$ on $M$
agree with the action of $A$ on $M \otimes A$).\par
Let $M$ be a partial right $A$--module (i.e. a partial abelian $A$--group).
Let us consider the ordinary tenzor product $M \otimes_{\bf Z} A$ of  abelian group $M$ by the ring  $A$ over {\bf Z}. Then the tensor $A$-completion of the partial $A$-module $M$ can be obtained by factorizing  $M \otimes_{\bf Z} A$ by all relations of the type
$(x \alpha, \beta) - (x, \beta)\alpha, x \in M, \alpha, \beta \in A.$ \par

 The following proposition is a copy of the corresponding result from the theory of modules. 
\begin{pr}
\label{pr3}
\begin{enumerate}
\item Let $M$ be a torsion--free abelian group, and $A$ a ring with 
a torsion--free additive group. Then $i: M \longrightarrow 
M \otimes A$ is injective, and $M \otimes A$ is torsion--free.
\item Let $M$ be a partial $A$--module without elements of order
2, and $A$ a ring whose additive group has no elements of order 2. 
Then $M \otimes A$ has no elements of order 2.
\end{enumerate}
\end{pr}
\par\medskip

{\it Non-abelian case}
\par\medskip
Before getting into details for the noncommutative case, we will briefly
describe the plan of an attack. By the definition of an $A$--group, the centralizers of
$G^A$ must be $A$--modules. This gives us a hint. The main idea is
to extend all centralizers of $G$ up to $A$--modules, adding a minimal
number of new relations. Centralizers in a CSA--group are abelian,
and we already know their $A$--completions.
But we need to have this tensor completion of given maximal abelian
subgroups $M$ of $G$ agree with all the other elements. The best
way to do that is to extend the centralizers of $G$ by the 
$A$--module $M \otimes A$  under the canonical monomorphism
$i: M \longrightarrow M \otimes A$. Then, using conjugate separation
of the subgroup $M$, we can extend the action of $A$ on $M \otimes A$
to all conjugates of $M$ (in order to satisfy axiom 3 of an $A$--group,
according to which the action of $A$ on $G$ commutes with all 
conjugations). To make this plan work, we need to put some natural
restrictions on $G$ and $A$, namely, to assume that the homomorphism
$i: M \longrightarrow M \otimes A$ is injective, and $A^+$
contains no elements of order 2 (i.e. $A^+$ is a CSA$^\ast$--group).

Before proving the main theorem, let us describe the construction of 
{\it complete tensor 
extension of centralizers of $G$ by the ring $A$.} We will assume
that $G$ is a non-abelian CSA--group and all maximal abelian subgroups
of $G$ are faithful over the ring $A$. \par
Let ${\cal C}=\{C_i, \ i \in I \}$ be the set of all centralizers 
of $G$ which  
satisfies the following conditions (S):
\begin{itemize}
\item any $C \in {\cal C}$ is not an $A$--module
(as a partial abelian $A$--subgroup in $G$);
\item any two centralizers from ${\cal C}$ are not conjugate in $G$;
\item any centralizer in $G$ which is not an $A$--module is 
conjugate to an element of ${\cal C}$.
\end{itemize}
 Then we can form the set 
${\cal H}_A=\{C_i \otimes A \mid i\in I\}$ and the set of canonical embeddings 
$\Phi_A = \{ \phi_i : C_i \hookrightarrow  C_i \otimes A \mid i\in I\}$.
Consider the tree extension of centralizers of this special type:
$$G^\ast = G({\cal C}, {\cal H}_A, \Phi_A)=G({\cal C}, A)$$
Now we can iterate this construction up to level $\omega$:
\begin{equation}
G=G^{(0)} < G^{(1)} < G^{(2)} < \cdots < G^{(n)} < \cdots
\label{eq:urav}
\end{equation}
where $G^{(n+1)}=G^{(n)}({\cal C}_n, A)$, and the set ${\cal C}_n$
of centralizers in the group $G^{(n)}$ satisfies the condition (S).
\begin{opr}
The union $\bigcup_{n \in \omega}G^{(n)}$ of the chain (\ref{eq:urav})
is called {\it complete tensor extension of centralizers of $G$ by the
ring $A$}. 
\end{opr}
\begin{theorem}
Let $A \geq  {\bf Z}$ be a ring without additive elements of order 2, and let a nonabelian partial
$A$--group $G$ be a CSA--group. If all maximal abelian subgroups 
of $G$ are faithful over $A$, then the tensor completion of $G$ by $A$ is
the complete tensor extension of centralizers of $G$ by the ring $A$.
\end{theorem}
$\Box$ Let us prove by induction on $n$ that each group $G^{(n)}$
from the chain (\ref{eq:urav}) meets the conditions of the theorem
and satisfies some special universal property. So now we will focus
our attention on the tree extension of centralizers $G^\ast = G({\cal C},
A)$, which was described above.\par
By the construction of a tree extension (see section 2), $G^\ast$
is the union of the chain of groups
$$
G=G_0 \hookrightarrow G_1 \hookrightarrow \cdots \hookrightarrow
G_\alpha \hookrightarrow \cdots \hookrightarrow G_\delta = G^\ast
$$
where the set of centralizers ${\cal C}= \{C_\alpha \mid \alpha < \delta \}$
is well-ordered, $C_\alpha = C_G(v_\alpha)$, and by definition: \par
$G_{\alpha + 1} = G_\alpha(v_\alpha, C_\alpha \otimes A);$ \par
$G_\gamma = \bigcup_{\alpha < \gamma}G_\alpha$ for a limit ordinal $\gamma$.
\begin{lemma}
For any $\ 1 \neq g \in G$ the centralizer $C_{G^\ast}(g)$ of $g$ in $G^\ast$
is isomorphic to $C_G(g) \otimes A$. 
\end{lemma}

$\Box$ By the choice of the set ${\cal C}$ the centralizer 
$C_G(g)$ is conjugate in $G$ to some centralizer 
$C_\alpha = C_G(v_\alpha) \in {\cal C}$, or $C_G(g)$  is an $A$--module
in $G$. In the first case, the centralizer $C_G(v_\alpha)$ has to be
extended exactly once on the level $\alpha$, and by corollary 2 of
lemma 1, the extended centralizer is equal to $C_\alpha \otimes
A = C_{G^\ast}(v_\alpha)$. Therefore the centralizer $C_{G^\ast}(g)$
in $G^\ast$ is conjugate to $C_{G^\ast}(v_\alpha)$, and hence
$C_{G^\ast}(g) \simeq C_G(g) \otimes A$. In the second case
$C_G(g)$ is an $A$--module, so $C_{G^\ast}(g) \simeq C_G(g) \otimes A$.
\hfill $\Box$

\begin{lemma}
\label{le1}
$G^\ast$ is a partial $A$--group, and the action of $A$ on $G$ is defined 
in $G^\ast$.
\end{lemma}

$\Box$ By lemma 4, in the group $G^\ast$ the centralizer $M=C_{G^\ast}(g)$
of any nontrivial element $g \in G$ is isomorphic to the group
$C_G(g) \otimes A$ which has a canonial structure of an $A$--module.
Using this isomorphism, the action of $A$ on $C_G(g) \otimes A$ can be
induced on $M$. By theorem 7, $G^\ast$ is a CSA--group, hence
$M$ is conjugately separated, and we can extend this action correctly
to all conjugates of $M$ in $G^\ast$ to satisfy axiom 3) from the definition
of an $A$--group. All maximal abelian subgroups in $G^\ast$ are
separated, so the intersections $M_1 \cap M_2$ are trivial, and
we have a well--defined action of $A$ on the group $G$ and all its 
conjugates. Moreover, axioms 1)--4) from the definition of an $A$--group hold
in $G^\ast$. \hfill $\Box$
\begin{lemma} 
\label{le2}
$G^\ast$ is a CSA--group, and all its maximal abelian subgroups are
faithful over $A$.
\end{lemma}

$\Box$ From the general theory of free products with amalgamation 
(see \cite{mks}) we know that elements of finite order lie in conjugates
of factors. By proposition~\ref{pr3}, the tensor completion $M \otimes 
A$ of any maximal abelian subgroup $M$ of $G$ has no elements of order 2.
Therefore by transfinite induction on the construction of the tree
extension $G^\ast$ we conclude that $G^\ast$ has no elements of order 2. 
$G^\ast$ is a CSA--group by theorem 7. Moreover, each abelian
subgroup of $G^\ast$ is either an $A$--module, or infinite cyclic.
In both cases, they are faithful over $A$. \hfill $\Box$

\begin{lemma}
The group $G^\ast$ has the following universal property with respect to the
canonical embedding $i: G \hookrightarrow G^\ast$: for any $A$--group
$H$ and any partial $A$--homomorphism $f: G \longrightarrow H$ 
there exists a partial $A$--homomorphism $f^\ast: G^\ast \longrightarrow
H$ such that $f=f^\ast \circ i$.
\end{lemma}

$\Box$ Let $M$ be a maximal abelian subgroup of $G$, and $M \in {\cal C}$.
Then $f(M)$ generates an abelian $A$--subgroup 
$N$ in $H$. Using the universal property of tensor completion for
abelian groups, we can find a homomorphism $\psi_M: M \otimes A 
\longrightarrow N$ 
such that $f=\psi_M \circ i_M$, where $i_M: M \hookrightarrow M \otimes A$
is the canonical embedding. 
By proposition 8, there exists $f^\ast : G^\ast \longrightarrow H$ with the
property $f=f^\ast \circ i$. \hfill $\Box$\par
\medskip
{\it Proof of the theorem.} Let $G^{(n+1)} = (G^{(n)})^\ast, \ G_0=G$,
where $G^\ast$ is defined as above. Then the union $\bigcup_{n \in \omega}
G^{(n)} = G \otimes A$
is the tensor $A$--completion of $G$. Indeed,
$G\otimes A$ is a partial $A$--group as a union of partial $A$--groups, but an
action of $A$ on $G \otimes A$ is defined everywhere: if $x \in G \otimes A$,
then
$x \in G^{(n)}$, and hence by lemma 5 the action of $A$ on $x$
is defined in $G^{(n+1)}$. So $G \otimes A$ is an $A$--group. \par
To prove that $G \otimes A$ is a tensor completion of $G$, we need to
verify the corresponding universal property.  Let $H$ be an
$A$--group and $f: G \longrightarrow H$ a homomorphism. Using lemma 7, 
we can extend $f$ to an $A$--homomorphism $f_n^\ast: G^{(n)} 
\longrightarrow H$ of partial $A$--groups. But then $f^\ast = \bigcup_{n
< \omega} f_n$ is an $A$--homomorphism such that $f=f^\ast \circ i$,
where $i$ is the canonical embedding $i:G \longrightarrow G \otimes A$.
\hfill $\Box$

\begin{co}
Let $G$ be a non-abelian CSA--group and $A \geq {\bf Z}$ a ring without additive elements of order 2. Then $G$ is faithful over $A$ iff each its abelian subgroup
is faithful over $A$.
\end{co}

\begin{co}
Let $G$ be an $A$--faithful non-abelian CSA--group and $A$ a ring without
additive elements of order 2. Then the tensor $A$--completion of $G$ is also
a CSA--group.
\end{co}

\begin{theorem}
Let $G$ be a torsion--free CSA--group and $A$ a ring with a torsion--free
additive group $A^+$. Then $G^A$ is a torsion--free CSA--group, and the 
canonical map $i: G \longrightarrow G^A$ is injective. In particular,
this is true for a hyperbolic torsion--free group $G$, a group
$G$ acting freely on a $\Lambda$--tree, and a universally free group
$G$.
\end{theorem}
$\Box$ A corollary of the previous theorem and proposition 12. 
\hfill $\Box$
\par

\section{A--free groups}

In this section we investigate some basic properties of $A$--free 
groups. First of all we establish some properties of $F^A$ as a CSA--group,
then introduce canonical forms of elements of $F^A$ and describe
commuting elements and conjugate elements in $F^A$.\par
Let $A$ be a ring with a torsion--free additive group $A^+$, and 
$F_A(X)$ an $A$--free group with base $X$. From theorem 2 we know
that $F_A(X)$ is the $A$--completion of the ordinary free group $F(X)=F$.
The group $F$ is a CSA--group, and all maximal abelian subgroups
of $F$ are faithful over $A$ (because they are infinite cyclic), therefore
the $A$--completion $F^A$ is the complete tensor extension of centralizers
of $F$ (theorem 8), and we can apply results of sections 4 and 5 to $F^A$.
In the case $A={\bf Q}$ G.Baumslag obtained similar results for the group
$F^{\bf Q}.$
\begin{theorem}
Let $F_A(X)$ be a non--abelian $A$--free group. Then 
\begin{enumerate}
\item the canonical mapping $F \longrightarrow F_A(X)$ is injective;
\item $F_A(X)$ is a torsion--free group with conjugately separated
maximal abelian subgroups;
\item $F_A(X)$ has no nonabelian solvable subgroups;
\item $F_A(X)$ has no nonabelian Baumslag -- Solitar subgroups;
\item $F_A(X)$ has trivial center and is directly undecomposable
\end{enumerate}
\end{theorem}

$\Box$ From theorem 8 and proposition 9. \hfill $\Box$
\par \medskip
By theorem 8, the group $F_A(X)$ is the union of the chain of subgroups:
$$
F=F^{(0)} < F^{(1)} < \cdots < F^{(n)} < \cdots
$$
where $F^{(n+1)} = F^{(n)}(C_n, A)$, and the set $C_n$ of centralizers
in $F^{(n)}$ satisfies the condition (S), i.e. any centralizer from
$C_n$ is infinite cyclic, no two of them are conjugate, and any
cyclic centralizer in $F^{(n)}$ is conjugate to one from $C_n$.
If $C(v)= \langle v \rangle \in C_n $ , then the element $v$ is called
{\it a root element of level n}; $V_n$ is the set of all root elements of level
$n$ such that only one generator of the subgroup $\langle v \rangle$
belongs to $V_n$; $V = \cup_{n \in \omega}V_n$. The {\it level function}
$\nu : F_A(X) \longrightarrow \omega$ is defined by the rule:
$$
\nu(g) = n \Longleftrightarrow g \in F^{(n)} \setminus F^{(n-1)},
$$
the number $\nu(g)$ is said to be {\it the level of the element g}. \par
Let us describe normal forms of elements of $F_A(X)$.\par
By section 4, if $\nu(g)=n$, then $g$ can be written in reduced form:
\begin{equation}
g = u_1v_1^{a_1}u_2v_2^{a_2} \cdots u_mv_m^{a_m}u_{m+1},
\label{reform}
\end{equation}
where $a_i \in A \setminus {\bf Z}$,
$v_i \in V_n$, $\nu(u_i) < n$, and if $v_i = v_{i+1}$, then 
$u_i \not \in v_i^{\bf Z}$. 
Moreover, we can assume that the elements $u_1, \ldots, u_{m+1}$
are also in reduced form as elements of $F^{(n-1)}$. This leads us
to an {\it explicit reduced form} of $g$. Similarly, following section 4,
one can introduce {\it canonical} and {\it semicanonical forms} of
elements and their {\it explicit} analogs.

\begin{pr} \label{pr20}
Let $g$ and $h$ be in reduced form (\ref{reform}):\par
$g = u_1v_{i_1}^{a_1} \cdots u_mv_{i_m}^{a_m}u_{m+1}, \ \ h = w_1v_{j_1}^{b_1} \cdots w_kv_{i_k}^{b_k}w_{k+1}.$ \par
Then $g=h$ if and only if 
\begin{enumerate}
\item[a)] $m=k$ and $i_s = j_s$, $a_s = b_s+n_s$
for some $n_s \in {\bf Z}$, $s = 1, \ldots, k$;
\item[b)] $u_{m+1} = v_{i_m}^{t_m}w_{m+1}, \ u_m = v_{i_{m-1}}^{t_{m-1}}
w_mv_{i_m}^{-t_m-n_m}, \ldots, u_1 = w_1v_{i_1}^{-t_1-n_1}$, where
$n_s$ is as in 1),$t_s \in {\bf Z}, \ s = 1, \ldots, m$.  

\end{enumerate}
\end{pr}

\begin{pr}
Suppose that $g$ and $h$ are in reduced form as in proposition~\ref{pr20}
and $a_s=b_s$ for all $s$. Then $g=h$ iff 
$g$ can be obtained from $h$ by a finite number of commutations
of the type
$v^av^t = v^tv^a, \ a \in A, \ v \in V_n, \ t \in {\bf Z}.$ 
\end{pr}
The following proposition gives us a description of commuting elements
in $F_A(X)$ which is completely analogous to the case of ordinary free 
groups.

\begin{pr}Let $F_A(X)$ be a nonabelian $A$--free group. Then:
\begin{enumerate}
\item elements $x, \ y \in F_A(X)$ commute if $x=z^a, \ y=z^b$
for some $z \in F_A(X)$ and $a,b \in A$;
\item centralizers in $F_A(X)$ are exactly free $A$--modules of rank 1.
\end{enumerate}
\end{pr}

$\Box$ By theorem 8, $F_A(X)$ is the complete tensor extension of 
centralizers by the ring $A$. According to this construction, if
$g \in F_A(X)$ is an element of level $\alpha$, then by lemma 4
$C(g) \simeq {\bf Z} \otimes_{\bf Z} A \simeq A$, and the centralizer
$C(g)$ is not extended on the next levels. Hence 2) is proved, and
1) follows from 2). \hfill $\Box$\par \medskip

The reduced form (\ref{reform}) of an element $g$ is called
{\it cyclically reduced} if any cyclic permutation of this form
is also reduced.

\begin{theorem}
Let $g$ be a cyclically reduced element from the $A$--free group
$F_A(X)$. Then 
\begin{enumerate}
\item if $g$ is conjugate to an element $h$, then $h$ has the same level as $g$, 
and the level of the conjugating element is not greater than that of $g$;
\item if $g$ is conjugate to a cyclically reduced element $p_1 \cdots p_n$,
$n \geq 2$, of level $k$, then $g$ can be obtained by permuting
$p_1, \ldots, p_n$ cyclically and then conjugating by a root
element of level $k-1$.
\end{enumerate}
\end{theorem}

$\Box$ Follows from lemma 2 from section 6 and the construction 
of the complete tensor extensions of centralizers of $F$. \hfill $\Box$

\section{Open problems}

\begin{question}
Using extensions of centralizers, construct the $A$--completions
of free solvable groups and groups of the type $F/[N,N]$
(these groups are not CSA--groups).
\end{question}

\begin{question}
Prove that the $A$--completion of a free nonabelian nilpotent
(solvable) group can be non--nilpotent (non--solvable).
\end{question}

The affirmative answer to question 3 will show the necessity to
define the notion of $A$--completion relatively to a variety of groups.
\begin{question}
Describe the tensor $A$--completion relatively to the varieties of 
nilpotent and solvable groups.
\end{question}
Let us mention that Hall's $A$--completion of a nilpotent group could be
different from that involved in question 4. \par
\medskip
Finally we formulate some algorithmic problems for free $A$--groups over
an effective \cite{rab} integral domain $A$ of characteristic 0.
\begin{question}
Prove that in the $A$--free group $F^A$
the word problem and the conjugacy problem are algorithmically
decidable.
\end{question}

\begin{question}
Is the Diophantine problem over an $A$--free group decidable?
\end{question}

\begin{question}
Is the universal theory of an $A$--free group decidable?
\end{question}
If the domain $A$ is residually {\bf Z}, then the answres to
questions 5, 6 and 7 are affirmative.

\end{document}